\documentclass[
11pt,%
final,%
tightenlines,%
notitlepage,%
twoside,%
onecolumn,%
nobibnotes,%
nofootinbib,%
superscriptaddress,%
noshowpacs,%
noshowydk,%
showdoi,%
centertags,%
maik]%
{revtex4n}

\input jpack_e
\setcaptionmargin{5mm}
\onelinecaptionstrue

\columntovsizetrue
\allowdisplaybreaks

\begin{document}

\newtheorem{teo}{Theorem}
\newtheorem*{teon}{Theorem}
\newtheorem{lem}{Lemma}
\newtheorem*{lemn}{Lemma}
\newtheorem{prp}{Proposition}
\newtheorem*{prpn}{Proposition}
\newtheorem{ass}{Assertion}
\newtheorem*{assn}{Assertion}
\newtheorem{assum}{Assumption}
\newtheorem*{assumn}{Assumption}
\newtheorem{stat}{Statement}
\newtheorem*{statn}{Statement}
\newtheorem{cor}{Corollary}
\newtheorem*{corn}{Corollary}
\newtheorem{hyp}{Hypothesis}
\newtheorem*{hypn}{Hypothesis}
\newtheorem{con}{Conjecture}
\newtheorem*{conn}{Conjecture}
\newtheorem{dfn}{Definition}
\newtheorem*{dfnn}{Definition}
\newtheorem{problem}{Problem}
\newtheorem*{problemn}{Problem}
\newtheorem{notat}{Notation}
\newtheorem*{notatn}{Notation}
\newtheorem{quest}{Question}
\newtheorem*{questn}{Question}

\theorembodyfont{\rm}
\newtheorem{rem}{Remark}
\newtheorem*{remn}{Remark}
\newtheorem{exa}{Example}
\newtheorem*{exan}{Example}
\newtheorem{cas}{Case}
\newtheorem*{casn}{Case}
\newtheorem{claim}{Claim}
\newtheorem*{claimn}{Claim}
\newtheorem{com}{Comment}
\newtheorem*{comn}{Comment}

\theoremheaderfont{\it}
\theorembodyfont{\rm}

\newtheorem{proof}{Proof}
\newtheorem*{proofn}{Proof}

\selectlanguage{english}
\Rubrika{\relax}
\CRubrika{\relax}
\SubRubrika{\relax}
\CSubRubrika{\relax}
%

\def\JournalNumber{0}
\def\JournalVolume{00}
%
%
%
\nameVolumeRus{}
\CnameVolumeRus{}
\nameIssueRus{\No}
\CnameIssueRus{}
\namePartRus{}
\namePagesRus{}
\nameYearShortRus{}
\JournalNameRus{}
\TranslitJournalNameRus{}
\JournalName{Regular and Chaotic Dynamics}
\JournalISSNCode{1560-3547}
\IssuePrice{}
\TransYearOfIssue{0000}
\TransCopyrightYear{2016}%
\OrigYearOfIssue{}
\OrigCopyrightYear{2016}%
\OrigIssueNo{\JournalNumber}
\OrigVolumeNo{\JournalVolume}
\TransVolumeNo{\JournalVolume}
\TransIssueNo{\JournalNumber}
\TransPartNo{}
\SHORTjournalPREFIX{RCD} 
\LONGjournalPREFIX{RegDyn} 
\BatFileName{call make_ps.bat} 
\BatSwitch{3} 
\IssueName{}
\SupplementNumber{}
\PublicationSerialNumberInYear{0}
\PublicationSerialNumberInVolume{0}
\ConditionalIssueDate{"year","month","day","name","type"}
\PagePrefix{}
\JournalISSNonlineCode{}
\JournalISSNCodeRus{}
\JournalISSNonlineCodeRus{}
\VolumeName{}
\IssnoName{none}
\PartnoName{}
\FpageNamepp{}
\FpageNnamep{}
\FpagePrefix{}
\LpageNnamepp{}
\LpageNamep{}
\LpagePrefix{}
\VolumePageNumbering{}
\JournalPubID{}
\FirstJournalPageNumber{}
\LastJournalPageNumber{}
\makeatletter
\def\MAIKlogo{RCD Editorial Office}
\def\maikpraefix{10.0000/S}
\edef\@ContentsHeadLineB{Simultaneous English language translation of the journal is available from \noexpand\MAIKlogo}
\def\Distributed{Distributed worldwide by Springer. }
\def\ArticlePages#1{\relax}
\@ifxundefined\CONT@sw{\@booleantrue\CONT@sw}{}%
\@booleantrue\showPACS@sw%
\@booleantrue\showKEYS@sw %
\@booleantrue\noOrigJournalVersion@sw
\@booleantrue\noOrigVolumeNo@sw
\@booleanfalse\noTransVolumeNo@sw
\makeatother
\input maikdoi %

\beginpaper


\input engnames
\titlerunning{Eight restricted problem}
\authorrunning{R. Lara and A. Bengochea}
\toctitle{A restricted four-body problem for the eight figure choreography}
\tocauthor{R. Lara and A. Bengochea}
\title{A restricted four-body problem for the eight figure choreography}
\firstaffiliation{
}%
\articleinenglish 
\PublishedInRussianNo
\author{\firstname{Ricardo}~\surname{Lara}}%
\email[E-mail: ]{lararicardo386@gmail.com}
\affiliation{
Department of Mathematics, ITAM\\
R\'io Hondo 1, Ciudad de M\'exico 01080, M\'exico}%
\author{\firstname{Abimael}~\surname{Bengochea}}%
\email[E-mail: ]{abimael.bengochea@itam.mx}
\begin{abstract}
In this work we introduce a planar restricted four-body problem where a massless particle moves under the gravitational influence due to three bodies following the eight figure choreography, and we explore some symmetric periodic orbits of this system which turns out to be non autonomous. We demonstrate that certain families of involutions are reversing symmetries of certain $N$--body problem. We consider a particular case of such families of reversing symmetries and we use it to study both theoretically and numerically certain type of symmetric periodic orbits of this system. The symmetric periodic orbits (initial conditions) were determined by means of solving some boundary value problems.
\end{abstract}
\received{October 13, 2019}
\accepted{Month XX, 20XX}%
\maketitle

\textmakefnmark{0}{)}%

\section{Introduction}
\label{I}

The $N$--body problem consists in the study of the dynamics of $N$ particles interacting through the Newton's law of gravity. Among the different type of solutions that appear in this system, the periodic orbits play a fundamental role since the system does not posses equilibrium solutions. In particular, a solution of the three-body, the eight figure solution \cite{moore,chenciner2000} gave rise to the study of a new type of solutions, the so called choreographies \cite{simo}. In a few words, a choreography is a periodic solution of the $N$--body problem where the bodies follow the same path with a common and constant time shift. One of the remarkable properties of the eight figure choreography is that is KAM stable, unlike other choreographies \cite{tomasz}. The eight figure choreography has been studied extensively in different perspectives, for instance curvature \cite{fujiwara}, or in the sense of bifurcations and stability \cite{galan,roberts}. The number of studies on choreographies has increased, up to now several choreographies are known for $n>3$ bodies \cite{kapela,simo}, and they have been studied in spaces of constant curvature \cite{montanelli1,montanelli2}, and for homogeneous potentials different from the gravitational \cite{hiroshi}. There are still open questions \cite{chenciner2003}, for instance, if choreographies with different masses could exist \cite{chenciner2004}.

In this article, we introduce a restricted four-body problem, related with the eight figure solution. We consider the study of the motion of a massless particle in presence of three bodies which follow the eight figure choreography. After setting such a restricted four-body problem, we focus on the periodic orbits of the system. In the classical restricted circular three and four body problems there exists a first integral, the so called Jacobi constant \cite{szebehely1964,burgos}, that allows the existence of monoparametric periodic solutions in the rotating frame, for a fixed value of the mass parameter. In our case the system does not have a first integral, therefore the periodic orbits are isolated, as it happens in the restricted elliptic three-body problem \cite{szebehely1964,broucke}. For the study of the periodic orbits we use certain reversing symmetries of the four body problem, which are compatible with the eight figure choreography. We prove the existence of a reversing symmetry for the planar $N$--body problem with $N = 2n+k$, for what is required that each pair of the $2n$ bodies have equal masses, without restriction in the masses of $k$ bodies. Later, we consider specifically the reversing symmetry for $n=1$, $k=2$, that will be useful for the restricted four-body problem. We compute numerically several periodic orbits by means of solving autonomous, and non-autonomous, boundary value problems.

As a matter of fact, our restricted four-body problem presents an inherent difficulty, different from those that appear in the circular and elliptic restricted three-body problems, which is that we do not have explicit formulas for the orbit of the primaries, we mean the eight figure choreography, although we know it through its initial condition. To the best of our knowledge, this restricted problem of four bodies has been not studied. We think that it could lead to the understanding of choreographies from a different point of view of those already considered.

The content of the article is the following. First, in section 2 we introduce the equations of motion of the $N$--body problem, and later, in Section 3 we state the definition of reversing symmetry, and a Theorem that gives sufficient conditions for the existence of periodic solutions in systems that possess reversing symmetries. In section 4 we give a proof of a reversing symmetry for the $N$--body problem where $N=2n+k$, in such a way each pair of the $2n$ bodies have equal masses, and in Section 4.1 we describe a specific case related with the restricted four-body problem. In Section 5 we define the restricted four-body problem, whereas in Section 6 we establish the equations that we study in order to determine numerically periodic orbits and we show the numerical results. In Section 7 we approximate the restricted four-body problem by one of two bodies, which can be considered if the test particle is far away from the primaries, and then we finish with the conclusions of this work.

\section{Equations of motion}
\label{II}

Consider $N$ particles on the plane with positive masses $m_i$, $i=1,\cdots,N$ interacting through the Newton's gravitational force. We denote respectively the position and velocity vectors by 
$$
{\bf r}_i =(x_i,y_i), \quad {\bf v}_i = \frac{d{\bf r}_i}{dt}, \quad i=1,\cdots,N.
$$
We assume an inertial frame of reference with origin at the center of mass of the $N$ bodies. The kinetic energy of the system and the gravitational potential are
$$
K = \frac{1}{2} \sum_{i=1}^N m_i {\bf v}_i^2, \quad U = - \sum_{i < j}^N \frac{Gm_im_j}{r_{ij}},
$$
where $G$ is the universal gravitational constant, and $r_{ij}=\vert {\bf r}_i - {\bf r}_j \vert$ the relative distance between the particles $i$ and $j$. The Lagrangian of the system is
$$
L = K-U,
$$
so the Euler-Lagrange equations
$$
\frac{d}{dt} \frac{\partial L}{\partial {\bf v}_i} - \frac{\partial L}{\partial {\bf r}_i}  = {\bf 0}, \quad i=1,\cdots,N
\label{2-1}
$$
give rise to the equations of motion
\begin{equation}
\frac{d^2{\bf r}_i}{dt^2} = \sum_{j \ne i}^N \frac{Gm_j({\bf r}_j-{\bf r}_i)}{r_{ij}^3}, \quad i=1,\cdots,N.
\label{2-2}
\end{equation}
The system (\ref{2-2}) can be rewritten as
\begin{equation}
\begin{array}{l}
\displaystyle \frac{d{\bf r}_i}{dt} = {\bf v}_i,  \\[6pt]
\displaystyle \frac{d{\bf v}_i}{dt} = {\bf a}_i,
\end{array}
\label{2-3}
\end{equation}
where
$$
{\bf a}_i = \sum_{j \ne i}^N \frac{Gm_j({\bf r}_j-{\bf r}_i)}{r_{ij}^3}, 
$$
for $i \in \{1,\cdots,N\}$. Considering 
$$
{\bf u} = ({\bf r}_1, \cdots, {\bf r}_N, {\bf v}_1, \cdots, {\bf v}_N),
$$
the equation (\ref{2-2}) can be written as a $4N$--dimensional first-order system
\begin{equation}
\frac{d{\bf u}}{dt} = {\bf F}({\bf u}),
\label{2-3}
\end{equation}
where
$$
(F_{2i-1},F_{2i})= {\bf v}_i, \quad (F_{2N+2i-1},F_{2N+2i}) = {\bf a}_i, \quad i=1,\cdots,N.
$$
The system (\ref{2-3}), equivalent to (\ref{2-2}), is not defined if at least two bodies are in the same position, i.e., when there is collision. Let
$$
\Delta_{ij} = \{ ({\bf r}_1,\cdots,{\bf r}_N) \in \mathbb{R}^{2N} \enskip \vert \enskip {\bf r}_i = {\bf r}_j \},
\quad \Delta = \bigcup_{i<j} \Delta_{ij},
$$
the configuration space of the system is $\Lambda = \mathbb{R}^{2N} \backslash \Delta$, and $\Omega = \Lambda \times \mathbb{R}^{2N}$ the phase space. We remark that the initial position and velocity vectors ${\bf r}_{i}(0) \in \Lambda$, ${\bf v}_{i}(0) \in  \mathbb{R}^{2N}$, $i=1,\cdots,N$, define a solution of the $N$--body problem, or in an equivalent way the vector ${\bf u}(0) \in \Omega$.

\section{Reversing symmetries and periodic orbits}
\label{III}
The reversing symmetries \cite{lamb,meiss} have been used successfully for studying the periodic orbits of problems described by differential equations \cite{galan,munoz0,munoz2,beng2013,beng2015}. In the following we give a brief review of some useful results.

Given (\ref{2-3}), we say that an involution $R: \Omega \longmapsto \Omega$ is a reversing symmetry if
$$
\frac{dR({\bf u})}{dt} = -{\bf F} \circ R{\bf (u)}
$$
holds. We denote the set of fixed points of $R$ as
$$
\textnormal{Fix}(R) = \{u \in \Omega \enskip \vert \enskip Ru = u \},
$$
and we say that these points are reversible configurations.

If $R$ is a reversing symmetry, we say that a solution ${\bf u}(t)$ of (\ref{2-3}), defined in its maximal domain, is $R$--symmetric if its orbit is invariant under $R$. We remark that an orbit is periodic and $R$--symmetric if and only if intersects two points of Fix$(R)$. Moreover, if $R$ and $\widehat R$ are reversing symmetries then we say that a solution ${\bf u}(t)$ is $(R,\widehat R)$--symmetric if there exist $t_0<t_1$ in such a way ${\bf u}(t_0) \in \textnormal{Fix}(R)$ and ${\bf u}(t_1) \in \textnormal{Fix}(\widehat R)$. In this case we also say that$[t_0, t_1]$ is the basic domain.

The next Theorem establish a main property between reversing symmetries and the periodicity of the solutions of differential equations \cite{munoz2,vander}.\\

\begin{teo}[Mu\~noz-Almaraz, Vanderbauwhede]
Let ${\bf u}(t)$ be a $(R, \widehat{R})$--sym\-met\-ric solution with basic domain $[0, T_0]$. Then ${\bf u}(t)$ is defined for all $t \in \mathbb{R}$, and we have for all $t \in \mathbb{R}$ and all $m \in Z$ that
$$
\begin{array}{c}
{\bf u}(-t) = R{\bf u}(t), \\[6pt]
{\bf u}(t) = \widehat{R}{\bf u}(2T_0 - t), \\[6pt]
{\bf u}(2mT_0 + t) = (\widehat{R} \circ R)^m{\bf u}(t). \\
\end{array}
$$
If there exists some $M \in \mathbb{N}$ such that $(\widehat{R} \circ R)^M = id$, the solution is $(R, \widehat{R})$--symmetric with period $T=2MT_0$.
\end{teo}

\section{Reversing symmetry $\Phi_\theta$ and its fixed points}
\label{reversing}

In this Section we demonstrate that the transformation $\Phi_\theta$, to be defined, is a reversing symmetry of the planar $N$--body problem with $N=2n+k$, $n,k,\in\mathbb{N}$. This reversing symmetry has the restriction that $2n$ bodies must have the same value in their masses at  pairs, we mean $m_{2i}=m_{2i-1}$ for $i=1,\cdots,{n}$. There is not a restriction about the masses for the other $k$ bodies\footnote{With the aim of simplify the proof, we exclude two cases: $k=0$, $n \ge 1$, and $n=0$, $k \ge 2$.}. Thereafter, we determine their fixed points. This type of fixed points was mentioned in the four body problem by Broucke \cite{broucke2004}, for the case of four equal masses.

First, we introduce the notation. We use column vectors, including nested vectors. For instance, if ${\bf x} \in \mathbb{R}^{2m}$, ${\bf x}_i \in \mathbb{R}^2$, $i=1,\cdots,m$, then we could write
$$
{\bf x} =
\left(
\begin{array}{c}
{\bf x}_1 \\
\vdots \\
{\bf x}_m
\end{array}
\right).
$$
Moreover, given any transformation $A:\mathbb{R}^2 \to \mathbb{R}^2$, we use brackets $[\cdot]$ to indicate
$$
A[{\bf x}] = \left(
\begin{array}{c}
A{\bf x}_1 \\
\vdots \\
A{\bf x}_m
\end{array}
\right).
$$
Besides the usual state vectors ${\bf r}_i$, ${\bf v}_i$, $i=1,\cdots,N$,
$$
{\bf r} = \left(
\begin{array}{c}
{\bf r}_1 \\
\vdots \\
{\bf r}_N
\end{array}
\right), \enskip
{\bf v} = \left(
\begin{array}{c}
{\bf v}_1 \\
\vdots \\
{\bf v}_N
\end{array}
\right), \enskip
{\bf u} = \left(
\begin{array}{c}
{\bf r} \\
{\bf v}
\end{array}
\right), 
$$
we also consider
\begin{equation}
\begin{array}{l}
{\bf q} = \left(
\begin{array}{c}
{\bf r}_1 \\
{\bf r}_2 \\
\vdots \\
{\bf r}_{2n-1}\\
{\bf r}_{2n}
\end{array}
\right), \enskip
\widetilde {\bf q} = \left(
\begin{array}{c}
{\bf r}_2 \\
{\bf r}_1 \\
\vdots \\
{\bf r}_{2n}\\
{\bf r}_{2n-1}
\end{array}
\right), \enskip
{\bf Q} = \left(
\begin{array}{c}
{\bf r}_{2n+1} \\
\vdots\\
{\bf r}_{2n+k}
\end{array}
\right),\\ [35pt]
{\bf p} = \left(
\begin{array}{c}
{\bf v}_1 \\
{\bf v}_2 \\
\vdots \\
{\bf v}_{2n-1}\\
{\bf v}_{2n}
\end{array}
\right), \enskip
\widetilde {\bf p} = \left(
\begin{array}{c}
{\bf v}_2 \\
{\bf v}_1 \\
\vdots \\
{\bf v}_{2n}\\
{\bf v}_{2n-1}
\end{array}
\right), \enskip
{\bf P} = \left(
\begin{array}{c}
{\bf v}_{2n+1} \\
\vdots\\
{\bf v}_{2n+k}
\end{array}
\right).
\end{array}
\label{eq-a1}
\end{equation}
In terms of (\ref{eq-a1}), the equations of motion $\dot{\bf u} = {\bf F}({\bf u})$ become
\begin{equation}
\begin{array}{c}
\left(
\begin{array}{c}
\dot{\bf q} \\
\dot{\bf Q} \\
\dot{\bf p} \\
\dot{\bf P}
\end{array}
\right) = 
\left(
\begin{array}{c}
{\bf p} \\
{\bf P} \\
{\bf f}({\bf q},{\bf Q}) \\
{\bf g}({\bf q},{\bf Q})
\end{array}
\right), \\[25pt]
{\bf f}({\bf q},{\bf Q}) = 
\left(
\begin{array}{c}
{\bf a}_1 ({\bf q},{\bf Q}) \\
{\bf a}_2 ({\bf q},{\bf Q}) \\
\vdots \\
{\bf a}_{2n-1} ({\bf q},{\bf Q}) \\
{\bf a}_{2n} ({\bf q},{\bf Q})
\end{array}
\right), \enskip
{\bf g}({\bf q},{\bf Q}) = 
\left(
\begin{array}{c}
{\bf a}_{2n+1} ({\bf q},{\bf Q}) \\
\vdots \\
{\bf a}_{2n+k} ({\bf q},{\bf Q})
\end{array}
\right),
\end{array}
\label{eq-a2}
\end{equation}
where 
\begin{equation}
{\bf a}_i ({\bf q},{\bf Q}) = \sum_{j=1, j \ne i}^{N} \frac{Gm_{j}({\bf r}_{j}-{\bf r}_{i})}{|{\bf r}_{j}-{\bf r}_{i}|^3}, \enskip i=1,\cdots,N.
\label{eq-a3}
\end{equation}
We also define rotation and reflection matrices in $\mathbb{R}^2$, namely
$$
G_{\theta} =
\left(
\begin{array}{cc}
\cos \theta & -\sin \theta \\
\sin \theta & \cos \theta
\end{array}
\right), \enskip
K =
\left(
\begin{array}{cc}
1 & 0 \\
0 & -1
\end{array}
\right),
\label{eq-a4}
$$
and the transformations $\widetilde P$, $\widetilde G_\theta$, $\widetilde K: \mathbb{R}^{4N} \to \mathbb{R}^{4N}$,
\begin{equation}
\begin{array}{c}
\widetilde P
\left(
\begin{array}{l}
{\bf q} \\
{\bf Q} \\
{\bf p} \\
{\bf P}
\end{array}
\right) = 
\left(
\begin{array}{l}
\widetilde {\bf q} \\
{\bf Q} \\
\widetilde {\bf p} \\
{\bf P}
\end{array}
\right), \enskip
\widetilde G_\theta
\left(
\begin{array}{l}
{\bf q} \\
{\bf Q} \\
{\bf p} \\
{\bf P}
\end{array}
\right) = 
\left(
\begin{array}{l}
G_\theta[{\bf q}] \\
G_\theta[{\bf Q}] \\
G_\theta[{\bf p}] \\
G_\theta[{\bf P}]
\end{array}
\right),\\[25pt] \relax
\widetilde K
\left(
\begin{array}{c}
{\bf q} \\
{\bf Q} \\
{\bf p} \\
{\bf P}
\end{array}
\right) = 
\left(
\begin{array}{c}
K[{\bf q}] \\
K[{\bf Q}] \\
-K[{\bf p}] \\
-K[{\bf P}]
\end{array}
\right).
\end{array}
\label{eq-a5}
\end{equation}
With this, we are able to define the reversing symmetry $\Phi_\theta$.\\

\begin{teo}
The transformation  $\Phi_\theta: \mathbb{R}^{4N} \to \mathbb{R}^{4N}$, $\Phi_\theta = \widetilde P \circ \widetilde G_\theta \circ \widetilde K$, is a reversing symmetry of (\ref{eq-a2}). 
\end{teo}

\proofn
First, we show that $\Phi_\theta$ is an involution. In accordance with this purpose, we give some properties. Given the identity $I_m$ in $\mathbb{R}^m$, we have for the inverse transformations
\begin{equation}
\begin{array}{c}
G_\theta G_{-\theta} = I_2, \enskip K^2 = I_2,\\[6pt]
\widetilde G_\theta \circ \widetilde G_{-\theta} = I_{4N}, \enskip \widetilde K \circ \widetilde K = I_{4N}, \enskip \widetilde P \circ \widetilde P = I_{4N}.
\end{array}
\label{eq-a6}
\end{equation}
In a similar way, we have the commutation relations
\begin{equation}
\begin{array}{c}
G_\theta K = K G_{-\theta},\\[6pt]
\widetilde G_\theta \circ \widetilde K = \widetilde K \circ \widetilde G_{-\theta}, \enskip \widetilde P \circ \widetilde K = \widetilde K \circ \widetilde P, \enskip \widetilde P \circ \widetilde G_\theta = \widetilde G_\theta \circ \widetilde P.
\end{array}
\label{eq-a7}
\end{equation}
By (\ref{eq-a5}) and (\ref{eq-a6}), it is straightforward to see that $\Phi_\theta \circ \Phi_\theta = I_{4N}$.

The second condition that $\Phi_\theta$ must fulfill is
\begin{equation}
\Phi_\theta \circ {\bf F} = - {\bf F} \circ \Phi_\theta
\label{eq-a8}
\end{equation}
For the left hand of (\ref{eq-a8}) we have
$$
(\Phi_\theta \circ {\bf F})
\left(
\begin{array}{c}
{\bf q} \\
{\bf Q} \\
{\bf p} \\
{\bf P}
\end{array}
\right) = \widetilde P
\left(
\begin{array}{c}
G_\theta K [{\bf p}] \\ \relax
G_\theta K [{\bf p}] \\ \relax
-G_\theta K [{\bf f}({\bf q},{\bf Q})] \\ \relax
-G_\theta K [{\bf g}({\bf q},{\bf Q})]
\end{array}
\right).
$$
Using $\widetilde P \circ \widetilde G_\theta \circ \widetilde K = \widetilde G_\theta \circ \widetilde K \circ \widetilde P$, the previous equation become
\begin{equation}
(\Phi_\theta \circ {\bf F})
\left(
\begin{array}{l}
{\bf q} \\
{\bf Q} \\
{\bf p} \\
{\bf P}
\end{array}
\right) =
\left(
\begin{array}{c}
G_\theta K [\widetilde {\bf p}] \\ \relax
G_\theta K [{\bf P}] \\ \relax
-G_\theta K [\widetilde {\bf f}({\bf q},{\bf Q})] \\ \relax
-G_\theta K [{\bf g}({\bf q},{\bf Q})]
\end{array}
\right),
\label{eq-a9}
\end{equation}
where
$$
\widetilde {\bf f}({\bf q},{\bf Q}) = 
\left(
\begin{array}{c}
{\bf a}_2({\bf q},{\bf Q}) \\
{\bf a}_1({\bf q},{\bf Q}) \\
\vdots \\
{\bf a}_{2n}({\bf q},{\bf Q}) \\
{\bf a}_{2n-1}({\bf q},{\bf Q})
\end{array}
\right).
$$
On the other hand, for the right side of (\ref{eq-a8}) we have
\begin{equation}
-({\bf F} \circ \Phi_\theta)
\left(
\begin{array}{c}
{\bf q} \\
{\bf Q} \\
{\bf p} \\
{\bf P}
\end{array}
\right) = 
\left(
\begin{array}{c}
G_\theta K [\widetilde {\bf p}] \\ \relax
G_\theta K [{\bf P}] \\ \relax
-{\bf f} (G_\theta K [\widetilde {\bf q}],G_\theta K [{\bf Q}]) \\ \relax
-{\bf g} (G_\theta K [\widetilde {\bf q}],G_\theta K [{\bf Q}])
\end{array}
\right).
\label{eq-a10}
\end{equation}
By comparing (\ref{eq-a9}) and (\ref{eq-a10}), we get that $\Phi_\theta$ is a reversing symmetry if
\begin{equation}
\begin{array}{l}
{\bf f} (G_\theta K [\widetilde {\bf q}],G_\theta K [{\bf Q}]) = G_\theta K [\widetilde {\bf f}({\bf q},{\bf Q})], \\ [6pt]
{\bf g} (G_\theta K [\widetilde {\bf q}],G_\theta K [{\bf Q}]) = G_\theta K [{\bf g}({\bf q},{\bf Q})],
\end{array}
\label{eq-a11}
\end{equation}
is satisfied. In order to show that (\ref{eq-a11}) holds, we rewrite the accelerations (\ref{eq-a3}); we distinguish three cases for the possible indices $1,\cdots,N$: odd and even indices for $1,\cdots,2n$, and $2n+1,\cdots,2n+k$. The odd case is given by
\begin{equation}
\begin{array}{l}
\displaystyle {\bf a}_{2i-1}({\bf q},{\bf Q}) = \sum_{j=1, j \ne i}^{n} \frac{G m_{2j-1} ({\bf r}_{2j-1}-{\bf r}_{2i-1})}{|{\bf r}_{2j-1}-{\bf r}_{2i-1}|^3}\\[10pt]
\displaystyle +\sum_{j=1}^{n} \frac{G m_{2j} ({\bf r}_{2j}-{\bf r}_{2i-1})}{|{\bf r}_{2j}-{\bf r}_{2i-1}|^3}
+\sum_{j=1}^{k} \frac{Gm_{2n+j}({\bf r}_{2n+j}-{\bf r}_{2i-1})}{|{\bf r}_{2n+j}-{\bf r}_{2i-1}|^3},
\end{array}
\label{eq-a12}
\end{equation}
with $i=1,\cdots,n$. For the even case we have
\begin{equation}
\begin{array}{l}
\displaystyle {\bf a}_{2i}({\bf q},{\bf Q}) = \sum_{j=1}^{n} \frac{G m_{2j-1} ({\bf r}_{2j-1}-{\bf r}_{2i})}{|{\bf r}_{2j-1}-{\bf r}_{2i}|^3}\\[10pt]
\displaystyle +\sum_{j=1, j \ne i}^{n} \frac{G  m_{2j}  ({\bf r}_{2j}-{\bf r}_{2i})}{|{\bf r}_{2j}-{\bf r}_{2i}|^3}
+\sum_{j=1}^{k} \frac{Gm_{2n+j}({\bf r}_{2n+j}-{\bf r}_{2i})}{|{\bf r}_{2n+j}-{\bf r}_{2i}|^3},
\end{array}
\label{eq-a13}
\end{equation}
where $i=1,\cdots,n$. For the other case, we have for $i=1,\cdots,k$,
\begin{equation}
\begin{array}{l}
\displaystyle {\bf a}_{2n+i}({\bf q},{\bf Q}) = \sum_{j=1}^{n} \frac{G m_{2j-1} ({\bf r}_{2j-1}-{\bf r}_{2i})}{|{\bf r}_{2j-1}-{\bf r}_{2n+i}|^3}\\[10pt]
\displaystyle +\sum_{j=1}^{n} \frac{G m_{2j}  ({\bf r}_{2j}-{\bf r}_{2n+i})}{|{\bf r}_{2j}-{\bf r}_{2n+i}|^3}
+\sum_{j=1, j \ne i}^{k} \frac{Gm_{2n+j}({\bf r}_{2n+j}-{\bf r}_{2n+i})}{|{\bf r}_{2n+j}-{\bf r}_{2n+i}|^3}.
\end{array}
\label{eq-a14}
\end{equation}
It is straightforward to see that ${\bf a}_i(\widetilde {\bf q},{\bf Q})$, for $i \in \{1,\cdots, N\}$ fixed, is obtained from ${\bf a}_i(q,Q)$ by means of using $m_{2j}=m_{2j-1}$, and interchanging, inside the summations, ${\bf r}_{2j-1} \leftrightarrow {\bf r}_{2j}$, for $j=1,\dots,n$. Doing this for (\ref{eq-a12}-\ref{eq-a14}), we have
\begin{equation}
\begin{array}{l}
{\bf a}_{2i}(\widetilde {\bf q},{\bf Q}) = {\bf a}_{2i-1}({\bf q},{\bf Q}), \enskip i=1,\cdots,n, \\[6pt]
{\bf a}_{2i-1}(\widetilde {\bf q},{\bf Q}) = {\bf a}_{2i}({\bf q},{\bf Q}), \enskip i=1,\cdots,n, \\[6pt]
{\bf a}_{2n+i}(\widetilde {\bf q},{\bf Q}) = {\bf a}_{2n+i}({\bf q},{\bf Q}), \enskip i=1,\cdots,k.
\end{array}
\label{eq-a15}
\end{equation}
Thus, by using (\ref{eq-a15}) in the first equation of (\ref{eq-a11}), we get
$$
{\bf f} (G_\theta K [\widetilde {\bf q}],G_\theta K [{\bf Q}]) = 
\left(
\begin{array}{c}
{\bf a}_1 (G_\theta K [\widetilde {\bf q}],G_\theta K [{\bf Q}]) \\
{\bf a}_2 (G_\theta K [\widetilde {\bf q}],G_\theta K [{\bf Q}]) \\
\vdots \\
{\bf a}_{2n-1} (G_\theta K [\widetilde {\bf q}],G_\theta K[{\bf Q}]) \\
{\bf a}_{2n} (G_\theta K [\widetilde {\bf q}],G_\theta K[{\bf Q}])
\end{array}
\right) \\[10pt] 
$$
$$
 =
\left(
\begin{array}{c}
G_\theta K [{\bf a}_1 (\widetilde {\bf q},{\bf Q})] \\ \relax
G_\theta K [{\bf a}_2 (\widetilde {\bf q},{\bf Q})] \\ \relax
\vdots \\ \relax
G_\theta K [{\bf a}_{2n-1} (\widetilde {\bf q},{\bf Q})] \\ \relax
G_\theta K [{\bf a}_{2n} (\widetilde {\bf q},{\bf Q})]
\end{array}
\right) \\[6pt]
$$
$$
=  
\left(
\begin{array}{c}
G_\theta K [ {\bf a}_2 ({\bf q},{\bf Q})] \\ \relax
G_\theta K [ {\bf a}_1 ({\bf q},{\bf Q})] \\ \relax
\vdots \\ \relax
G_\theta K [ {\bf a}_{2n} ({\bf q},{\bf Q})] \\ \relax
G_\theta K [ {\bf a}_{2n-1} ({\bf q},{\bf Q})]
\end{array}
\right) \\[6pt]
$$
$$
= G_\theta K [\widetilde {\bf f}({\bf q},{\bf Q})].
$$
The demonstration for the other equation of (\ref{eq-a11}) is similar.\qed

Once we have demonstrated that $\Phi_\theta$ is a reversing symmetry, in the following we determine the fixed points of $\Phi_\theta$. The equation that defines the set of fixed points of $\Phi_\theta$, we mean Fix$(\Phi_\theta)$, is
$$
(\widetilde P \circ \widetilde G_\theta \circ \widetilde K) ({\bf u}) = {\bf u}.
$$
In the previous equation a parameter (angle) appears. It is enough to determine the fixed points for $\Phi_0$ for characterize $\Phi_\theta$ since they are related by rotations, as we state in the following result.\\

\begin{prp}
Let ${\bf u}$ be a fixed point of $\Phi_\theta$, then $\widetilde G_\alpha ({\bf u}) $ is a fixed point of $\Phi_{2\alpha+\theta}$.
\end{prp}

\proofn
Let ${\bf u} \in \textrm{Fix}(\Phi_\theta)$, therefore 
$$
(\widetilde P \circ \widetilde G_\theta \circ \widetilde K) ({\bf u}) = {\bf u}.
$$
Multiplying both sides by $\widetilde G_\alpha$, and introducing the identity $I_{4N} = \widetilde G_{-\alpha} \circ \widetilde G_{\alpha}$ at the left of ${\bf u}$, the previous equation becomes
$$
(\widetilde G_\alpha \circ \widetilde P \circ \widetilde G_\theta \circ \widetilde K \circ \widetilde G_{-\alpha} \circ \widetilde G_\alpha) ({\bf u}) = \widetilde G_\alpha ({\bf u}),
$$
Taking into account the previous equation and (\ref{eq-a7}), we obtain
$$
(\widetilde P \circ \widetilde G_{2\alpha+\theta} \circ \widetilde K) (\widetilde G_\alpha ({\bf u})) = \widetilde G_\alpha ({\bf u}),
$$
therefore $\widetilde G_\alpha ({\bf u}) \in$ Fix$(\Phi_{2\alpha+\theta})$.\qed

According to the previous result, we can use $\textrm{Fix}(\Phi_0)$ for describing $\textrm{Fix}(\Phi_\theta)$. In order to determine $\textrm{Fix}(\Phi_0)$, it is convenient to consider vectors (\ref{eq-a1}). The set Fix$(\Phi_0)$ is conformed by the state vectors that satisfy the equation
\begin{equation}
K [\widetilde {\bf q}] = {\bf q}, \enskip K [{\bf Q}] = {\bf Q}, \enskip -K [\widetilde {\bf p}] = {\bf p}, \enskip -K [ {\bf P}] = {\bf P},
\label{eq-a16}
\end{equation}
which admits an infinite number of solutions. Considering ${\bf r}_{2i-1}$, ${\bf v}_{2i-1}$, $i=1,\cdots,n$, $x_{2n+l}$, $v_{y2n+l}$, $l=1,\cdots,k$, as independent parameters, the solution of (\ref{eq-a16}) is given by
\begin{equation}
\begin{array}{l}
x_{2i} = x_{2i-1}, \enskip y_{2i} = -y_{2i-1}, \enskip i=1,\cdots,n, \\[6pt]
v_{x2i} = -v_{x2i-1}, \enskip v_{y2i} = v_{y2i-1}, \enskip i=1,\cdots,n, \\[6pt]
y_{2n+i}=0, \enskip v_{x2n+i}=0, \enskip i=1,\cdots,k.
\end{array}
\label{revgeneral}
\end{equation}
Notice that bodies with indices $2i-1$ and $2i$, for $i \in \{1,\cdots,n \}$ fixed, appear at the vertices of an isosceles triangle, and its velocities are related by a reflection along the direction of the base of the triangle ($y$--axis). On the other hand, the bodies with index $2n+i$, $i \in \{1,\cdots,k \}$, lie along the $x$--axis, and their velocities are orthogonal to their positions.

\begin{rem}
If we permutate the indices in $\Phi_\theta$ for the bodies with same mass, we mean $2j \longleftrightarrow 2j-1$, for $j \in \{1,\cdots,n \}$, the new transformation is also a reversing symmetry. Something similar holds for Fix$(\Phi_\theta)$.
\end{rem}

\subsection{Reversible configurations for the restricted four-body problem}

In the following we present the reversible configurations that we have considered for studying specific orbits of the restricted four-body problem.

Let us consider the reversible configuration (\ref{revgeneral}) for $n=1$, $k=2$. We set that the bodies 1, 2 and 3 have unitary mass, and that the fourth be infinitesimal. As we mentioned, if we permutate the bodies 1 and 2 we will obtain a reversible configuration. Nevertheless, for our restricted four-body problem we can permutate three bodies, those with unitary mass. In our problem it is relevant to identify reversing symmetries related by a permutation, therefore from here on we add an index to the reversing symmetry described in Section 4. In our case only is useful the cyclic permutation $\sigma:1 \to 3 \to 2 \to 1$, so we consider three reversible configurations Fix$(\Phi_{0,j})$, $j=1,2,3$. For reference, we write down explicitly Fix$(\Phi_{0,1})$:
\begin{equation}
\begin{array}{l}
x_1 = x_2 , \enskip y_1 =  - y_2, \enskip  y_{3}=0, \enskip y_{4}=0, \\ [6pt]
v_{x1} = -v_{x2}, \enskip v_{y1} = v_{y2}, \enskip v_{x3}=0, \enskip v_{x4}=0.
\end{array}
\label{revt=0}
\end{equation}
The other two reversible configurations can be defined recursively. We have for $j=2,3$ that Fix$(\Phi_{0,j})$ is Fix$(\Phi_{0,j-1})$ subject to $\sigma$.

\section{Restricted problem for the eight figure}
\label{IV}
\begin{figure}[h]
\centering
\includegraphics[width=90mm]{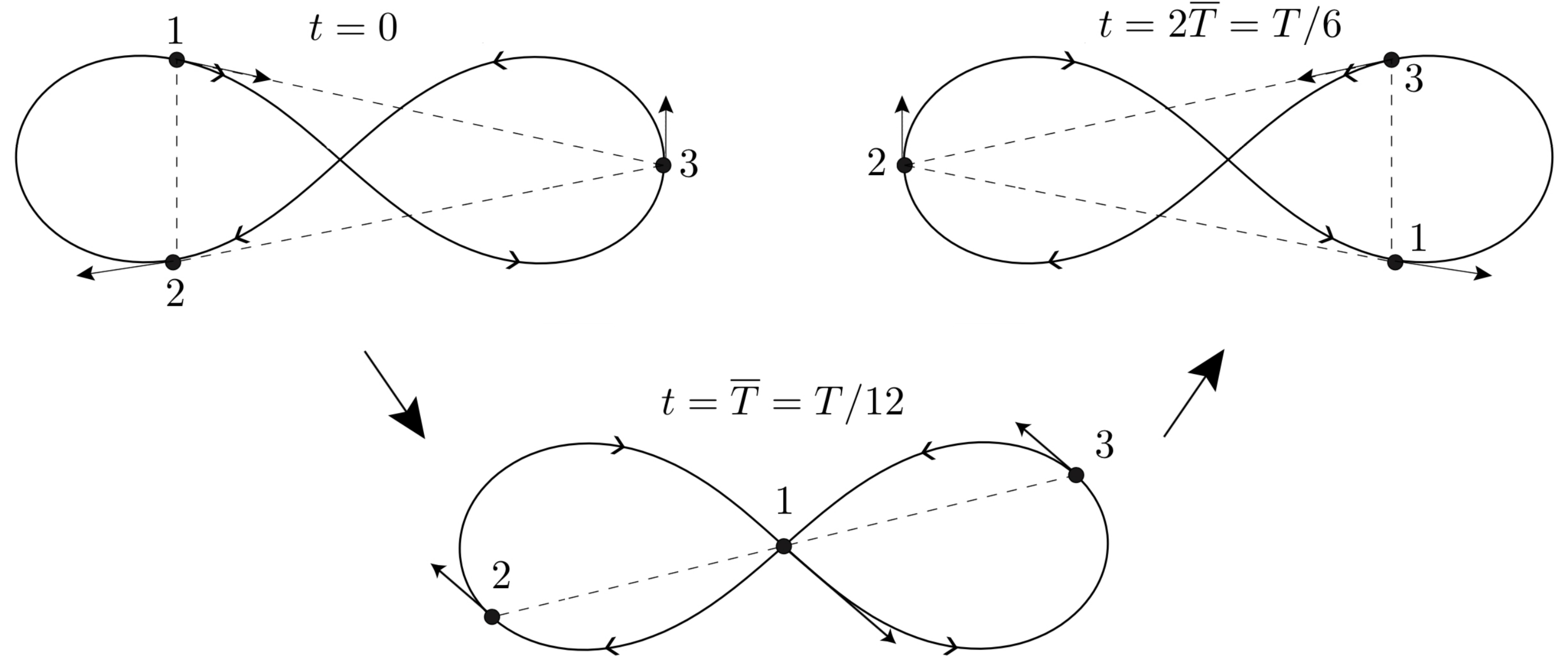}
\caption{Evolution of the eight figure choreography that starts at isosceles configuration.}
\label{fig1}
\end{figure}

Let us consider four particles in a fixed plane. Three of them, identified with indices $i=1,2,3$, follow the eight figure choreography, and a fourth one is massless moving by the influence of the others. We choose units \cite{chenciner2000} in such a way $m_i=1$, $i=1,2,3$, $G=1$, $T=12\overline{T} = 6.32591398$, where $T$ is the period of the choreography, and $2 \overline{T}$ the time that passes between two consecutive isosceles configurations (see Fig. ~\ref{fig1}).

We assume that at the initial time $t=0$, the particles $i=1,2,3$, appear in the first isosceles configuration depicted in Fig. ~\ref{fig1}, that is ${\bf r}_i(0) = {\bf r}_{ci}$, ${\bf v}_i(0) = {\bf v}_{ci}$, $i=1,2,3$, where
$$
\begin{array}{l}
{\bf r}_{c1} = (-0.54050854325, 0.3452633140), \\ [6pt]
{\bf r}_{c2} = (-0.54050854325, -0.3452633140), \\ [6pt]
{\bf r}_{c3} = (1.081017086500, 0), \\ [6pt]
{\bf v}_{c1} = (1.0971223818, -0.23360476285), \\ [6pt]
{\bf v}_{c2} = (-1.0971223818, -0.23360476285), \\ [6pt]
{\bf v}_{c3} = (0, 0.46720952570).
\end{array}
$$
The vectors ${\bf r}_i(t)$, ${\bf v}_i(t)$, $i=1,2,3$, are known functions (numerically) of time, $T$--periodic, and do not depend on the state vectors of the massless particle.

On the other hand, the motion of the fourth particle depends on the vectors ${\bf r}_i(t)$, $i=1,2,3$. The equations of motion of the fourth particle are
\begin{equation}
\ddot{\bf r}_4 = \sum_{j = 1}^3 \frac{{\bf r}_j(t)-{\bf r}_4}{r_{4j}^3}.
\label{3-1}
\end{equation}

It was already mentioned that one difficulty to study the equations (\ref{3-1}) is that we do not have an explicit expression for the position of the primaries, which does not happen in the classic cases, for instance the circular and elliptic restricted three-body problems.

\section{Periodic orbits in the restricted four-body problem}
\label{V}

With the objective of determining periodic orbits, we study those orbits which pass through fixed points of two reversing symmetries, possibly different. In our case the reversible configurations Fix$(\Phi_{0,j})$, $j=1,2,3$ (see Section 4.1) are compatible with the restricted four-body problem.

\begin{figure}[h]
\centering
\includegraphics[width=70mm]{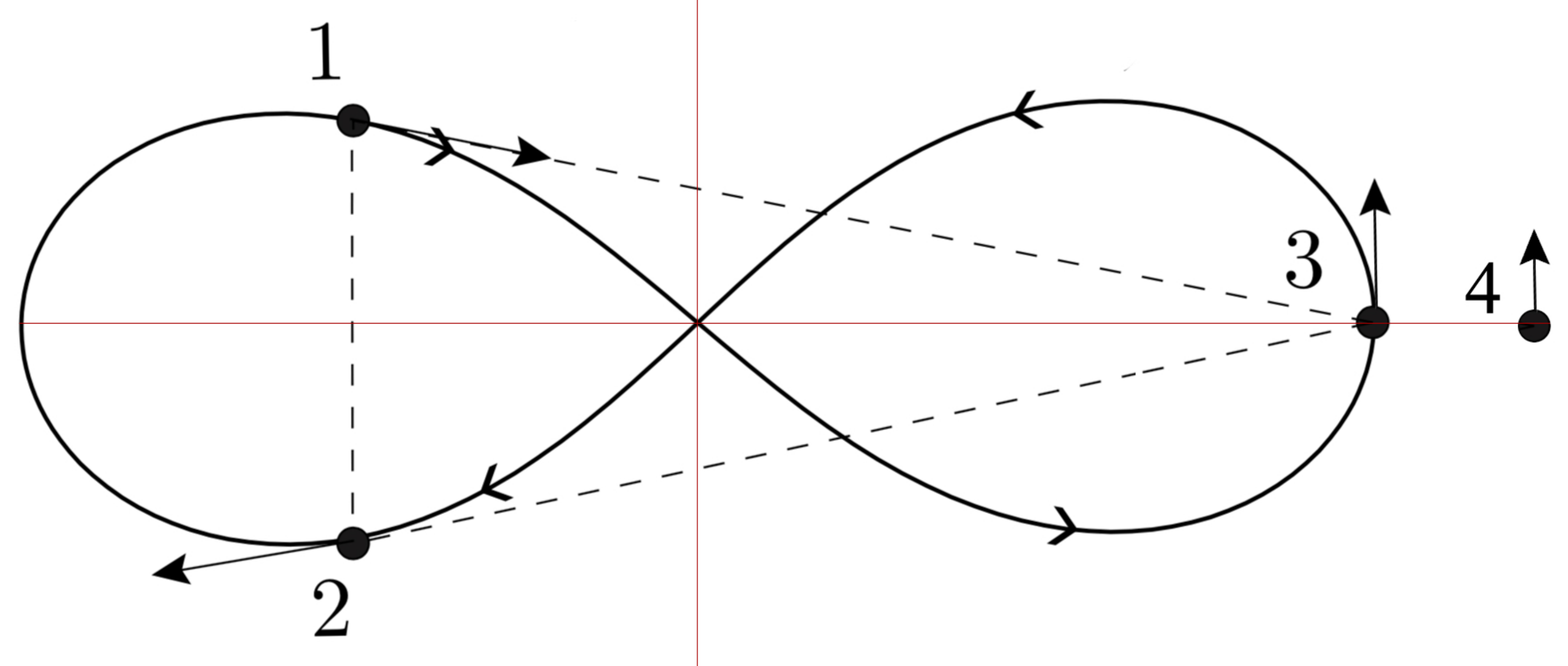}
\caption{Reversible configuration of the restricted four-body problem.}
\label{fig2}
\end{figure}

We state that at $t=0$ the orbits coincides with Fix$(\Phi_{0,1})$. Thus, the initial condition of the fourth particle matches Fix$(\Phi_{0,1})$, we mean, it is of the form $y_{40}=0$, $v_{x40}=0$, according to (\ref{revt=0}) (see Fig. \ref{fig2}). Now, we need that at some $T_0 \ne 0$ the orbit reaches a reversible configuration Fix$(\Phi_{0,j})$, $j \in \{1,2,3\}$. This only can happen if $T_0 = 2m\overline{T}$, $m \in \mathbb{N}_0$, because these are the times where the configuration of the three bodies is compatible with the eight figure choreography.

Therefore, we have to find unknowns $T_0 = 2m\overline{T}$, $m \in \mathbb{N}$, $x_{40},v_{y40} \in \mathbb{R}$, in such a way
\begin{equation}
\phi_{T_0}(x_{40},0,0,v_{y40}) = (x_{41},0,0,v_{y41}),
\label{num}
\end{equation}
where $\phi$ is the flow of (\ref{3-1}), and $x_{41},v_{y41}$, real numbers (its specific values are irrelevant, in the sense of periodicity). According to Theorem 1, if (\ref{num}) is satisfied then we have determined a periodic orbit since the reversible configurations Fix$(\Phi_{0,j})$, $j=1,2,3$ are related by a cyclic permutation in the indices of the primaries. 

The period of the orbit is determined by the value of $M$, according to Theorem 1. In order to identify $M$ we just need to see what reversing symmetries are associated to the orbit. By construction, at $t=0$ the orbit coincides with Fix$(\Phi_{0,1})$, therefore we only have three compositions for such reversing symmetries, we mean $(\Phi_{0,j} \circ \Phi_{0,1})$, $j=1,2,3$. It is straightforward to see that
$$
(\Phi_{0,1} \circ \Phi_{0,1}) = id, \enskip (\Phi_{0,2} \circ \Phi_{0,1})^3 = id, \enskip (\Phi_{0,3} \circ \Phi_{0,1})^3 = id.
$$
Thus, if the orbit passes through two fixed points of $\Phi_{0,1}$ then $M = 1$, and the corresponding period becomes $T = 2T_0$. On the other hand, if the orbit passes through a fixed point of $\Phi_{0,1}$ and $\Phi_{0,j}$ for some $j \in \{2,3\}$ then $T = 6T_0$.

With the aim of solve numerically this problem, we set boundary value problems (BVPs) associated to (\ref{3-1}) whose solutions lead to initial conditions which meet (\ref{num}). Our initial conditions are of the form $(x_{40},0,0,v_{y40})$ which for brevity we denote as $(x_{40},v_{y40})$.

\subsection{Boundary Value Problem - $T_0$ fixed}
\label{VA}

We establish two different BVPs, in each one of them $T_0$ is fixed as an multiple of $\overline{T}$. For these problems we we will have two unknowns, that is $x_{40}$ and $v_{y40}$. In the first BVP we set that at the time $t = 2p \overline{T}$, $p \in \mathbb{N}$, the fourth particle lies on the horizontal axis, we mean (\ref{3-1}) restricted to
\begin{equation}
\begin{array}{cc}
y_4(0) = 0, & y_4(2p \overline{T}) = 0, \\ [6pt]
v_{x4}(0) = 0.
\end{array}
\label{bvp-1}
\end{equation}
The solutions of (\ref{bvp-1}) are defined by the set
\begin{equation}
C_{(y,2p)} = \left \{ (x_{40},v_{y40}) \in \mathbb{R}^2\, | \, y_4(2p\overline{T}) = 0 \right \}.
\label{cy}
\end{equation}
In a similar way, at the time $t = 2q \overline{T}$ the fourth particle has zero horizontal velocity, thus
\begin{equation}
\begin{array}{c}
\begin{array}{cc}
y_4(0) = 0, & v_{x4}(2q \overline{T}) = 0. \\ [6pt]
v_{x4}(0) = 0.
\end{array}
\end{array}
\label{bvp-2}
\end{equation}
The solutions of (\ref{bvp-2}) are described by
\begin{equation}
C_{(v_x,2q)} = \left \{  (x_{40},v_{y40}) \in \mathbb{R}^2 \,  | \, v_{x4}(2q\overline{T}) = 0 \right \}.
\label{cvx}
\end{equation}
Notice that, in a generic sense, (\ref{cy}) and (\ref{cvx}) define curves in the plane $x_{40}v_{y40}$. A relevant property of these sets is that for $p=q$ the intersection is conformed by points that satisfy (\ref{num}), we mean, if $(x_{40},v_{y40}) \in C_{(y,2p)} \cap C_{(v_x,2p)}$, then $(x_{40},0,0,v_{y40})$ gives rise to a periodic orbit.

\subsection{Boundary Value Problem - $T_0$ variable}
\label{VB}

Next, we define a BVP in such a way the characteristic time $T_0 \in \mathbb{R}_{> 0}$ is unknown. We remark that for this BVP the time $T_0$ does not need to be an even multiple of $\overline{T}$. Thus, we consider (\ref{3-1}) with the condition 
\begin{equation}
\begin{array}{cc}
y_4(0) = 0, & y_4(T_0) = 0,  \\ [6pt]
v_{x4}(0) = 0, & v_{x4}(T_0) = 0. 
\end{array}
\label{bvp-3}
\end{equation}
We define the set of solutions of (\ref{bvp-3}) as
\begin{equation}
C_{R} = \left \{  (x_{40},v_{y40},T_0) \in \mathbb{R}^2 \times  \mathbb{R}_{> 0} \, | \, y_4(T_0) = 0, \, v_{x4}(T_0) = 0 \right \}.
\label{cR}
\end{equation}
The set (\ref{cR}) is conformed by curves in the three-dimensional space $x_{40}v_{y40}T_0$. The points within $C_{R}$ such that $T_0 = 2m\overline{T}$, $m \in \mathbb{N}$, define symmetric periodic orbits because the corresponding solutions pass through reversible configurations at times $t=0,T_0$.

\subsection{Numerical results}
\label{VC}

In order to compute solutions of the BVPs previously defined, we need a starting solution of the corresponding problem, that we will call a ``seed''. We obtain the ``seeds'' by means of solving root finding problems, in the following we explain for $C_{(y,2p)}$. According to (\ref{cy}) we have $T_0 = 2p \overline{T}$, for some $p \in \mathbb{N}$, and we want to determine unknowns $x_{40},v_{y40}$, so that
\begin{equation}
y_4(T_0,x_{40},v_{y40}) = 0.
\label{root}
\end{equation}
We need an approximate value of the root $(x_{40},v_{y40})$, namely $(x_{4*},v_{y4*})$. If we fix one of the variables, for instance $x_{40} = x_{4*}$, we only need to find $v_{y40}$. If $(x_{4*},v_{y4*})$ is a good enough approximation, usually with a Newton's method, we will find a value $v_{y40}$ such that (\ref{root}) holds. Similar root problems can be established for $C_{(v_x,2q)}$ and $C_{R}$. Once we have a ``seed'', we are able to compute solutions of the corresponding BVP.

\begin{figure}[!h]
\centering
\includegraphics[width=85mm]{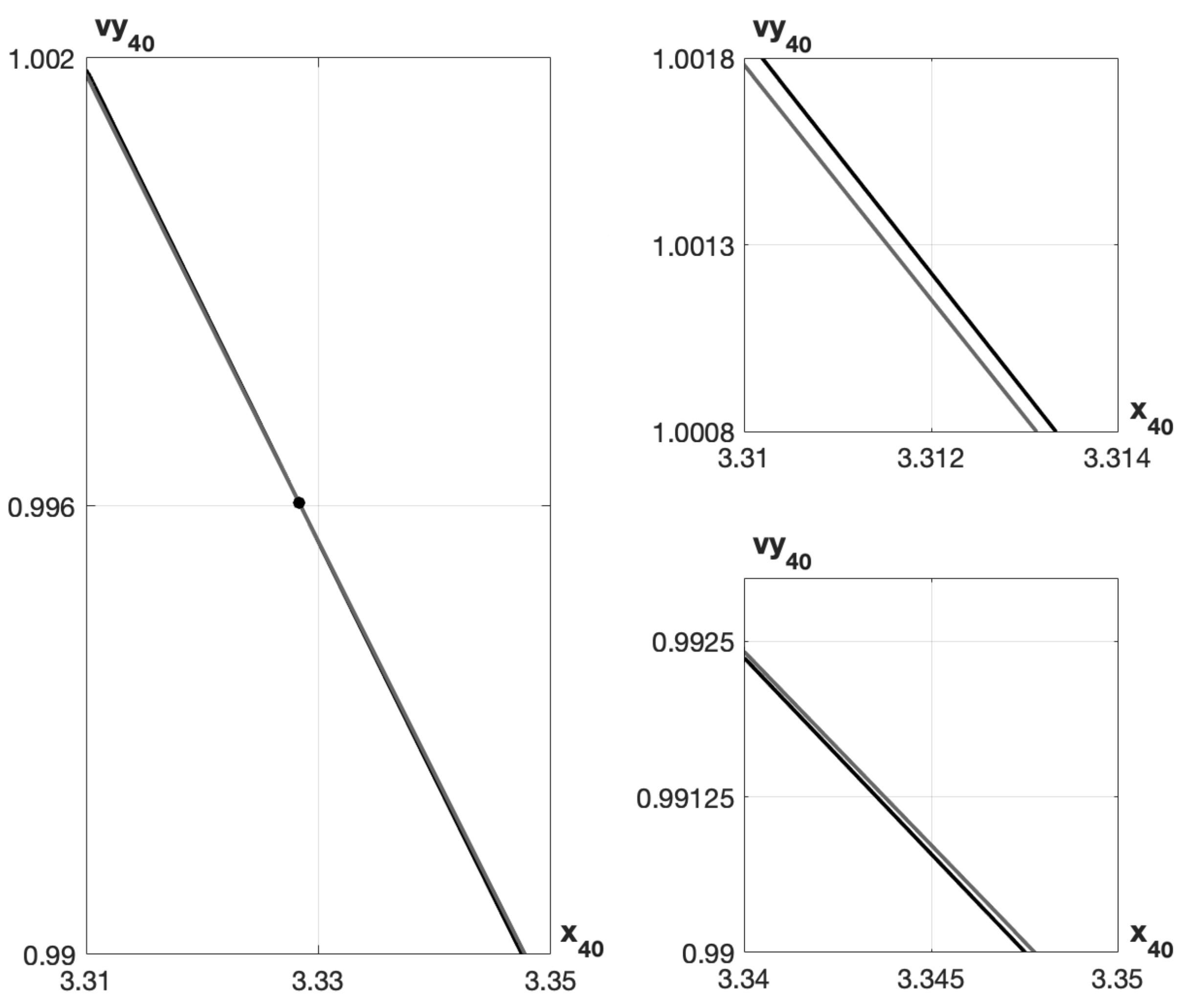}
\caption{Sets $C_{(y,20)}$, $C_{(v_x,20)}$ in black and gray series, respectively. At the intersection, indicated with a closed dot, a periodic orbit appears. The initial condition of this orbit is shown in Table 1, 15th row.}
\label{fig4}
\end{figure}

\begin{figure}[!h]
\centering
\includegraphics[width=85mm]{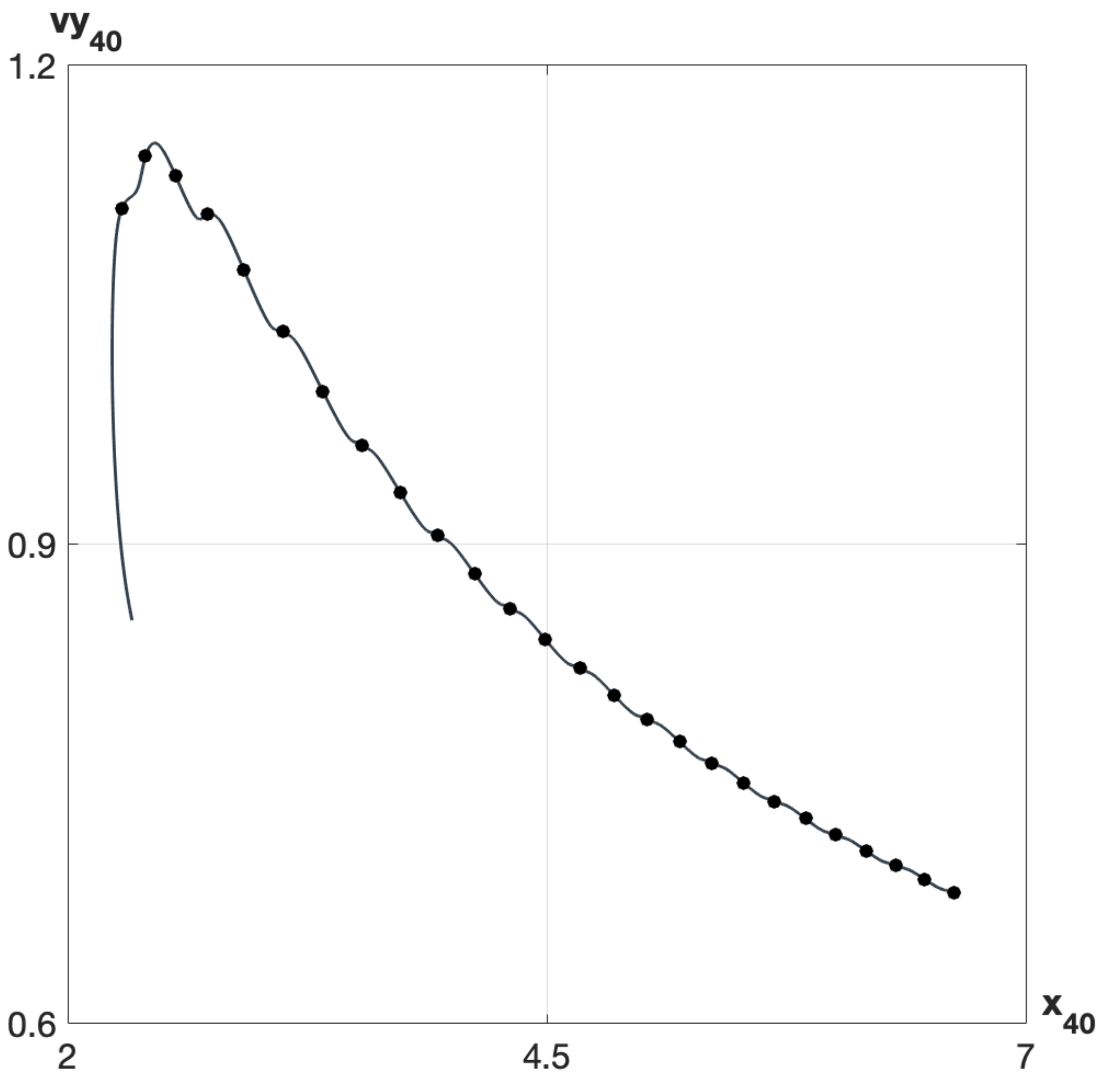}
\includegraphics[width=85mm]{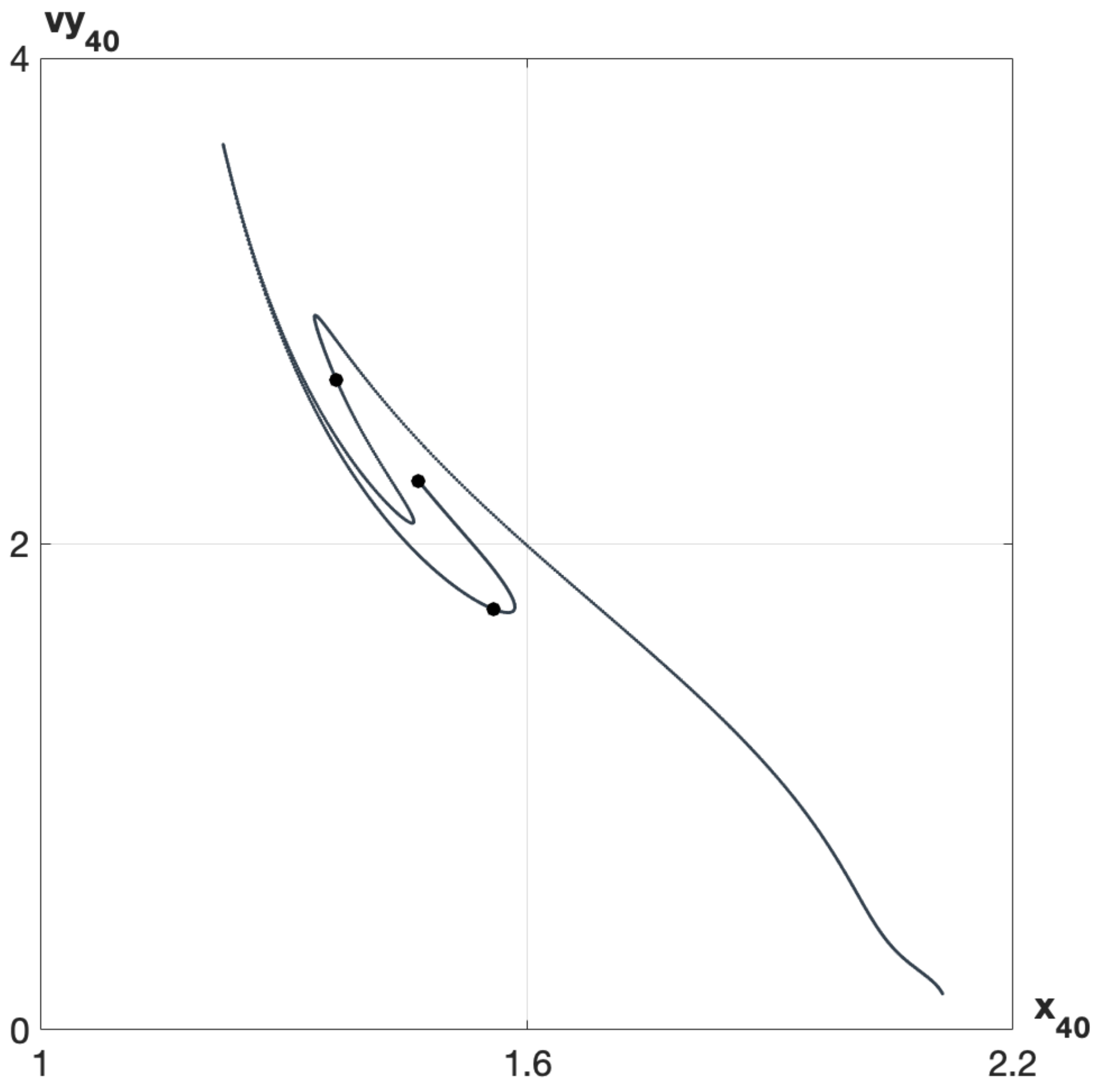}
\caption{Projection of the set $C_R$ on the plane $x_{40}v_{y40}$. The initial conditions represented by a closed dot define periodic orbits.}
\label{fig5}
\end{figure}

In Fig. ~\ref{fig4} we show a small part of the sets $C_{(y,20)}$, $C_{(vx,20)}$, and in Fig. ~\ref{fig5} we exhibit the projection of $C_{R}$ (computed part) in the plane $x_{40}v_{y40}$.

From the different BVPs, in the sense of periodicity, $C_R$ is the most practical because in this set we only need to identify the orbits in which $T_0$ is an even multiple of $\overline{T}$. In Fig. \ref{fig5}, first plot from top to bottom, we show one principal branch of $C_R$, which was computed without difficulties because the test particle does not pass close to the primaries. However, in the other plot within Fig. \ref{fig5} there are initial conditions that lead to orbits where the test particle passes near the primaries at some time. These orbits present numerical difficulties since the equations of motion are not defined at collisions.

\begin{figure}[!ht]
\centering
\includegraphics[width=72mm]{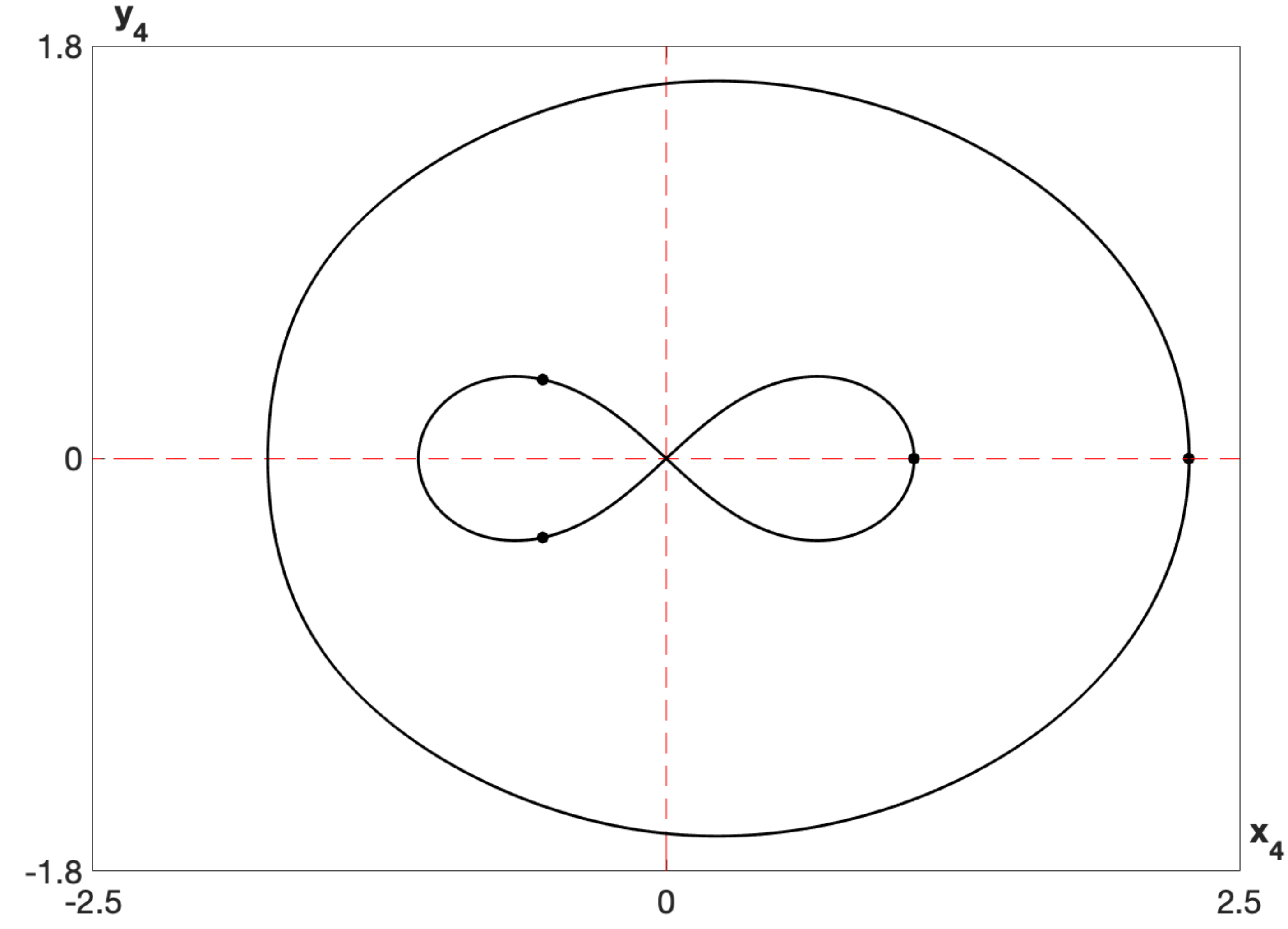} 
\includegraphics[width=78mm]{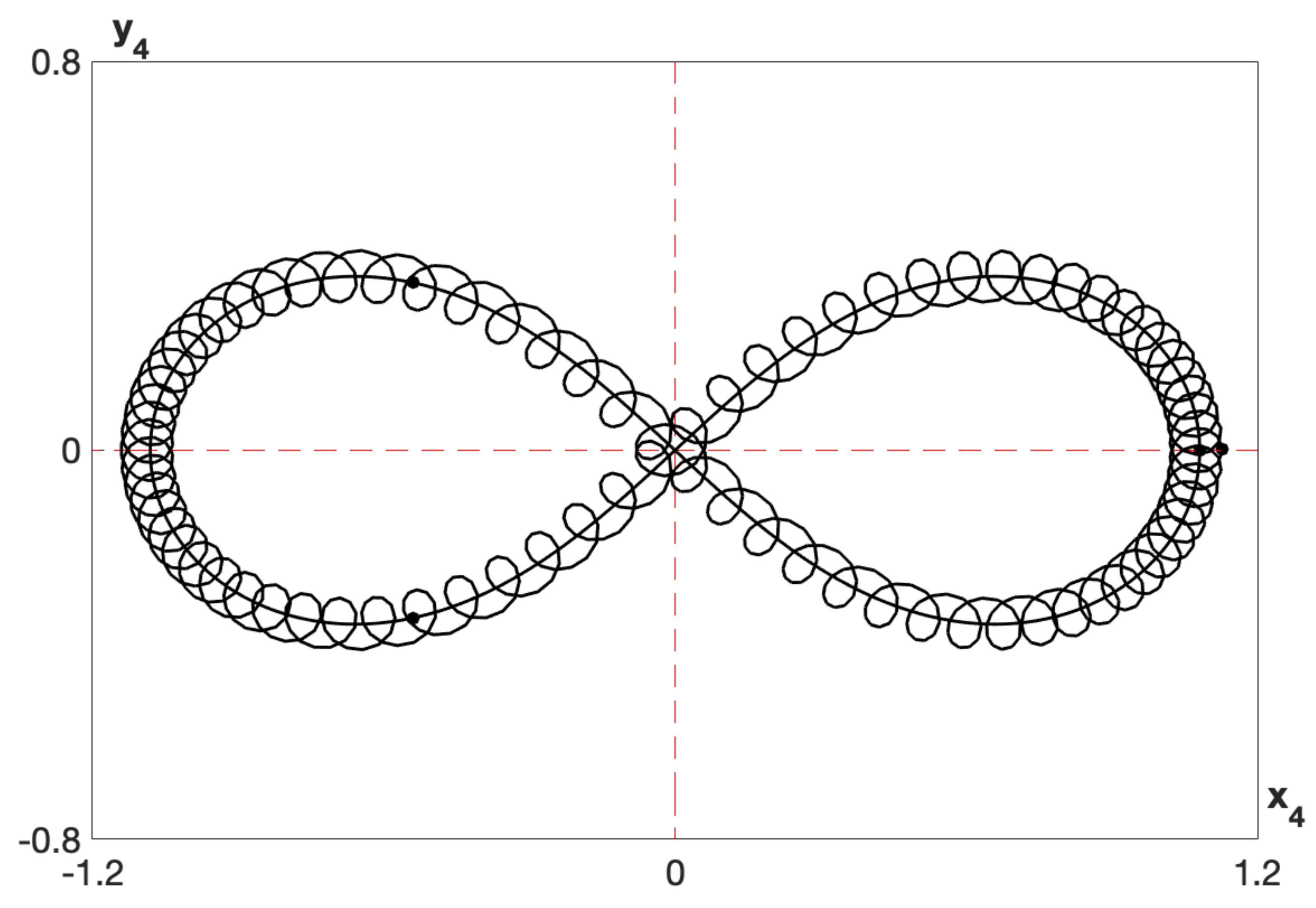}
\caption{From left to right, the orbits have period $48\overline{T}$ and $12\overline{T}$, respectively. The corresponding initial conditions are given in 9th and 3th rows of Table 1.}
\label{fig6}
\end{figure}

\begin{figure}[!ht]
\centering
\includegraphics[width=74mm]{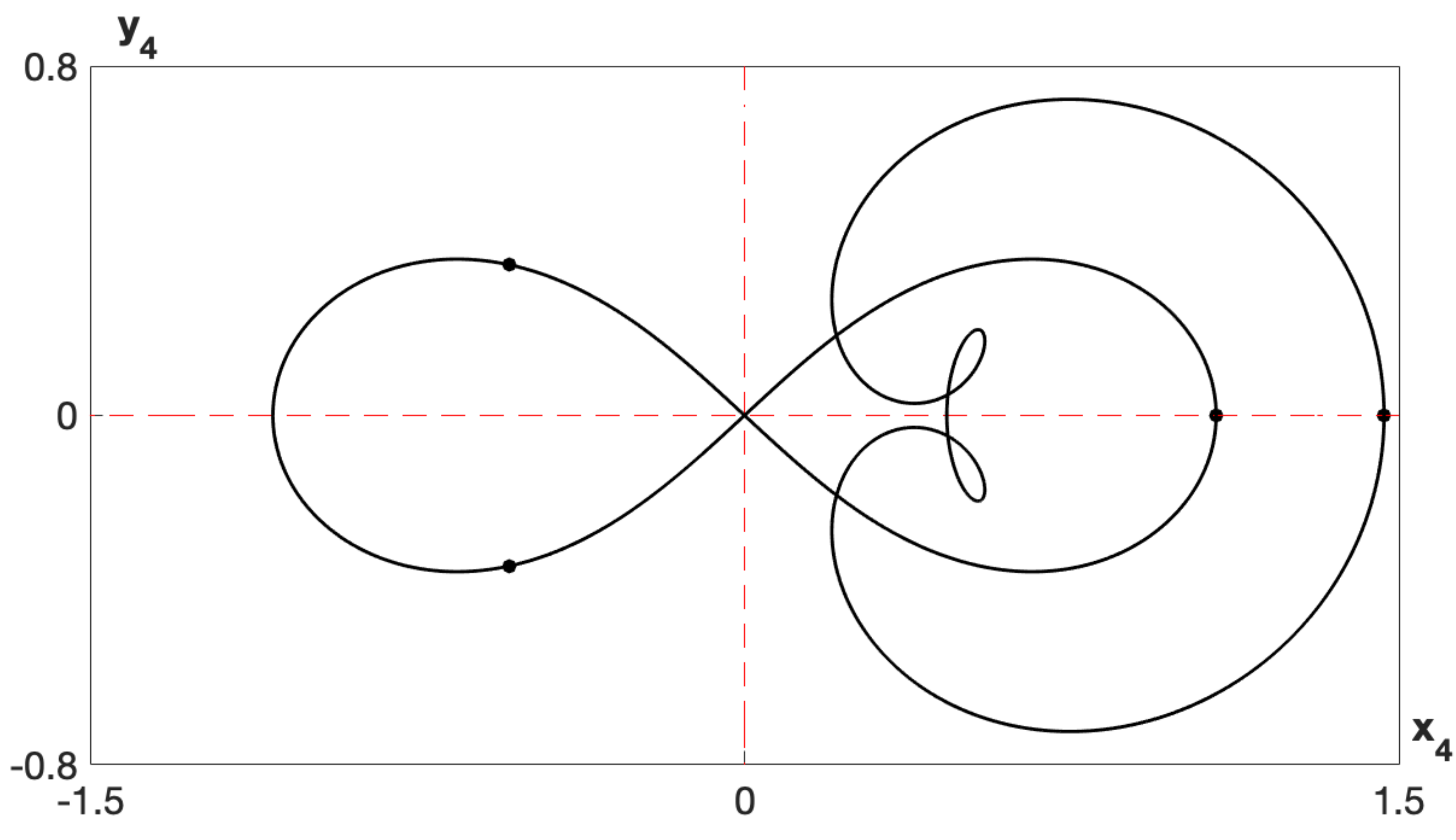} 
\includegraphics[width=74mm]{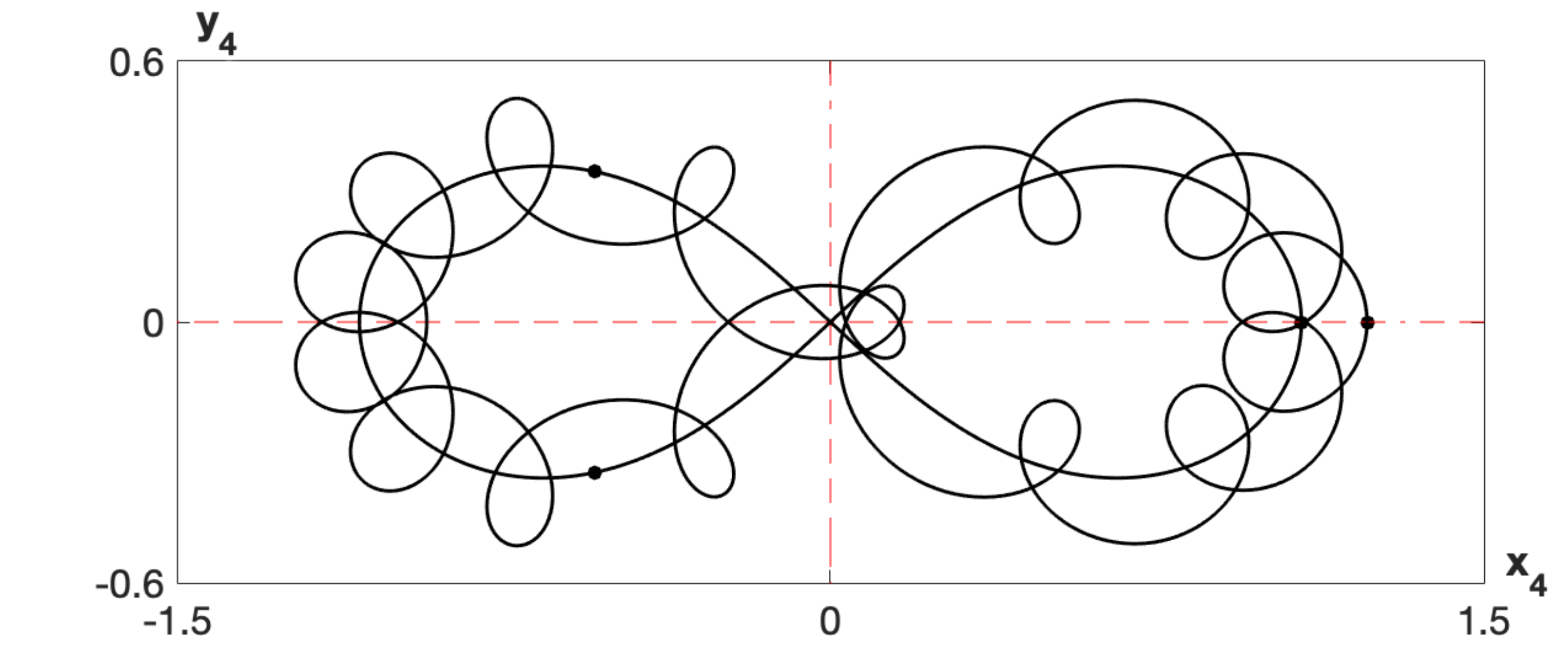} \\
\includegraphics[width=75mm]{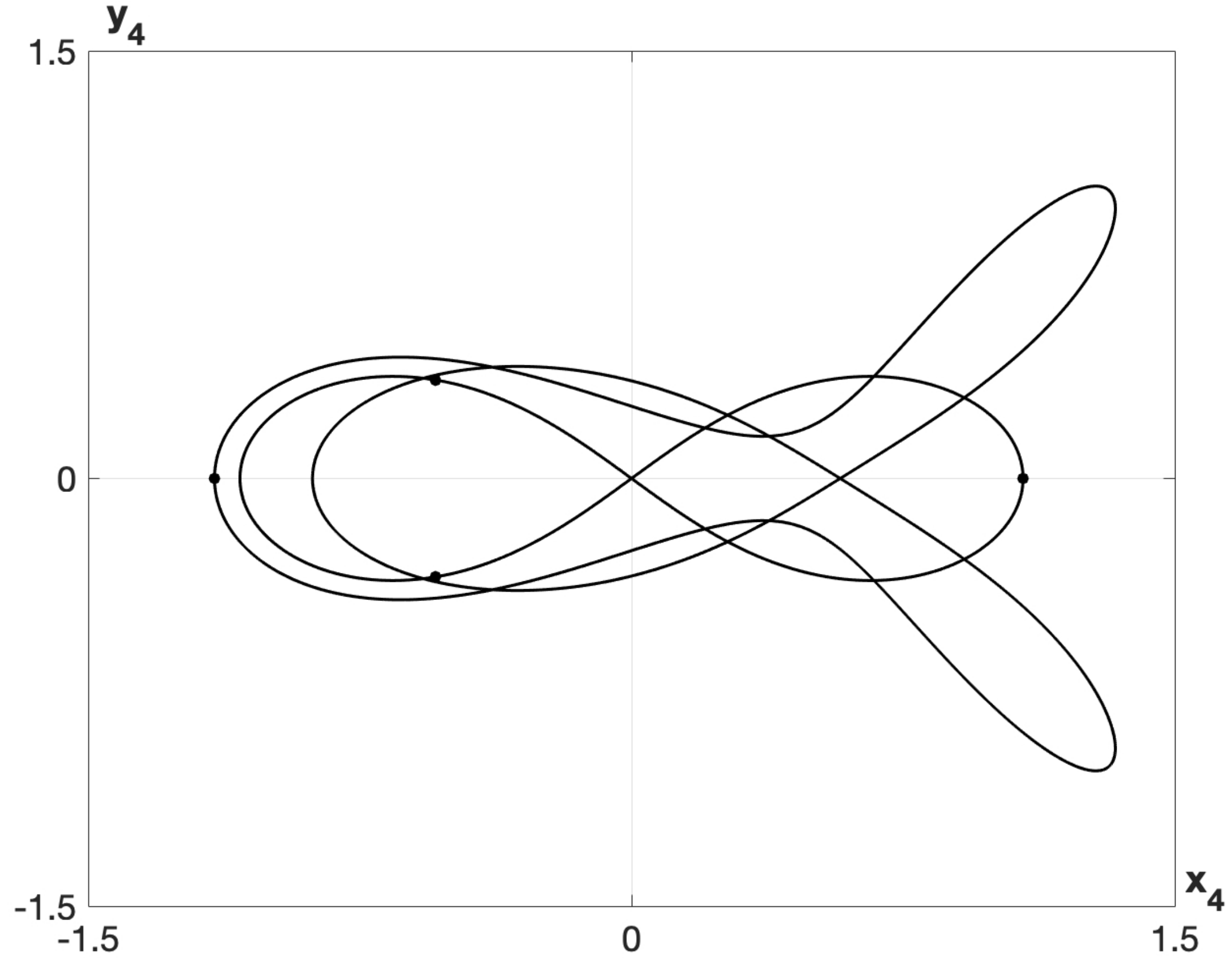} 
\includegraphics[width=78mm]{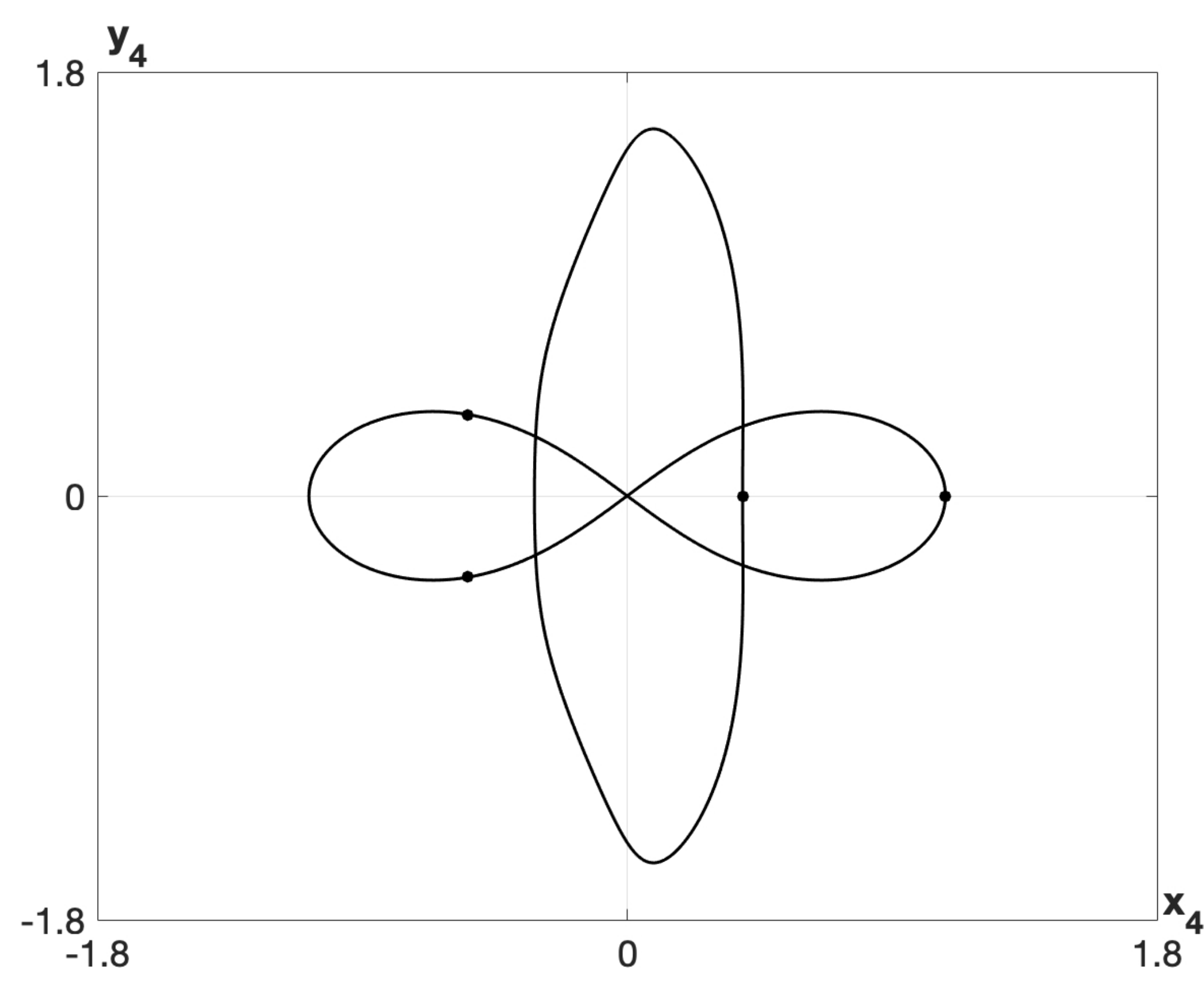} \\
\includegraphics[width=35mm]{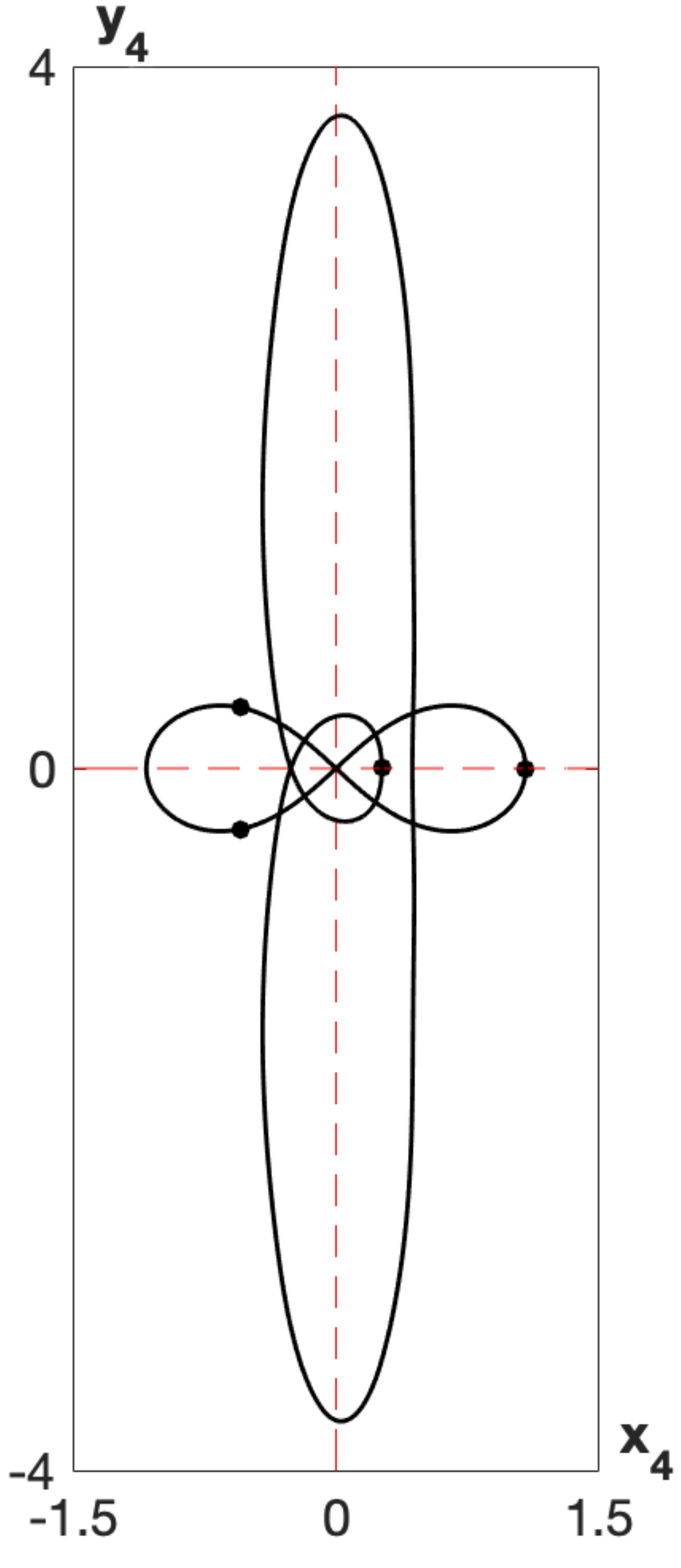}
\caption{From top to bottom, left to right, the first orbit has period $24\overline{T}$. The posterior orbits have periods $12\overline{T}$, $48\overline{T}$, $12\overline{T}$ and $120\overline{T}$, respectively. The corresponding initial conditions are given in 8th, 5th, 1th, 2th and 4th rows of Table 1.}
\label{fig7}
\end{figure}

We found symmetric periodic orbits with different properties. There are two kinds of orbits (limit cases) for those who the dynamics can be easily explained. For example, for certain small values of both $x_{40}$ and $v_{y40}$, the massless particle could turn around only one primary at all time, defining a binary system. On the other hand, for $x_{40}$ large enough, and adequate $v_{y40}$, the orbit of the fourth particle resembles an elliptic orbit (see next Section). Examples of these two limit cases are shown in Fig. \ref{fig6}. With the exception of the limit cases, in general the orbits present a complicated structure. For instance, there are orbits that in some part the test particle passes very close to the primaries, whereas in other part passes far away of the primaries. This can be appreciated in the orbits shown in Fig. \ref{fig7}. Notice that, as a consequence of the reversibility, all the periodic orbits computed are symmetric with respect to the horizontal axis.

\begin{table}[!ht]
\caption{Initial conditions for symmetric periodic orbits in the restricted four-body problem. The orbit passes through reversible configurations at times $t=0,T_0$. Here $12\widetilde{T}$ is the period of the eight figure choreography, and $T$ the period of the orbit of the restricted four-body problem.}
\begin{tabular}{|c|c|c|c|c|}
\hline
$\#$ & $T_0/\overline{T}$ & $T/\overline{T}$  & $x_4 $ & $vy_4$  \\
\hline
1 & 8 &  48 & -1.151372102323705 &  1.192735308310391 \\
2 & 6 &  12 & 0.392064354827005  & -2.088580677571261 \\
3 & 6 &  12 & 1.145057806500420  & 4.421084342099486 \\
4 & 20 & 120&0.261908739769502	& 0.768218486423285 \\
5 & 6 & 12& 1.231780839019731   & 3.098731109349930 \\
6 & 2 & 12& 1.364108936002170   & 2.676664885954700 \\
7 & 2 & 12& 1.557889835185201	& 1.732435901189350 \\
8 & 4 & 24& 1.465940005977230	& 2.259081336057390 \\
9 & 8 & 48& 2.280410388953660	& 1.110720371401547 \\
10 & 10 & 60& 2.403183107401021	& 1.143021410598030 \\
11 & 12 & 24& 2.559935679202217	& 1.131013972457270 \\
12 & 14 & 84& 2.727991976218170	& 1.106651149430582 \\
13 & 16 & 96& 2.917809831461645	& 1.071706543150445 \\
14 & 18 & 36& 3.120691699353664   & 1.033349295535902 \\
15 & 20 & 120&3.328354859295013	& 0.996039191967667 \\
16 & 22 & 132&3.533520671511424	& 0.962303086050867 \\
17 & 24 & 48& 3.734052676172197	& 0.932325325122930 \\
18 & 26 & 156&3.929756244211741	& 0.905635308058075 \\
19 & 28 & 168&4.120942435281265	& 0.881706024834630 \\
20 & 30 & 60& 4.307972049253366	& 0.860093902208141 \\
21 & 32 & 192&4.491174076273875	& 0.840442308932970 \\
22 & 34 & 204&4.670837162699625	& 0.822465134430957 \\
23 & 36 & 72& 4.847214688965527	& 0.805930897509638 \\
24 & 38 & 228&5.020530247022094	& 0.790650437478274 \\
25 & 40 & 240&5.190982479530962	& 0.776467608355671 \\
26 & 42 & 84& 5.358748764439068	& 0.763252315834217 \\
27 & 44 & 264&5.523988373388391	& 0.750895210415233 \\
28 & 46 & 276&5.686844913076257	& 0.739303645791943 \\
29 & 48 & 96& 5.847448514845465	& 0.728398528559842 \\
30 & 50 & 300&6.005917432292149	& 0.718111887323824 \\
31 & 52 & 312&6.162359623497286	& 0.708384923725778 \\
32 & 54 & 108&6.316873887698146	& 0.699166486221477 \\
33 & 56 & 336&6.469550993413046	& 0.690411825225105 \\
34 & 58 & 348&6.620474509569266	& 0.682081600922888 \\
\hline
\end{tabular}
\label{table1}
\end{table}

\section{Two-body approximation}
\label{VI}

This restricted four-body problem can be approximated by one of two bodies if the test particle is far enough away from those that follow the eight figure choreography. For the approximation we consider a fictitious body with mass $m=3$, position ${\bf r}$, and velocity ${\bf v}$, and the massless body ($m_4=0$) with position ${\bf r}_4$, and velocity ${\bf v}_4$. For these two bodies we consider an elliptic solution. From all the points that conform the ellipse we are interested exclusively in the apocenter and pericenter, since only at these configurations ${\bf r}_4$ and ${\bf v}_4$ fulfill the reversible configurations described in Section \ref{V}. Assuming the origin at the center of mass, so, at the apocenter the state vectors of the two-body problem are given by
$$
\begin{array}{c}
\displaystyle {\bf r} = -\frac{m_4}{m+m_4} a(1+e) \hat{\imath},  \enskip {\bf v} = - \frac{m_4}{\sqrt{m+m_4}} \sqrt{\frac{G(1-e)}{a(1+e)}} \hat{\jmath}, \\ [10pt]
\displaystyle {\bf r}_4 = \frac{m}{m+m_4} a(1+e) \hat{\imath},  \enskip {\bf v}_4 =  \frac{m}{\sqrt{m+m_4}} \sqrt{\frac{G(1-e)}{a(1+e)}} \hat{\jmath},
\end{array}
$$
where $a$ is the semi-major axis, $e$ the eccentricity, and $\hat{\imath}$, $\hat{\jmath}$ the usual canonical vectors. We are assuming that the origin of the coordinate system is located at the center of mass of the system, and that the semi-major axis goes in the direction $\hat{\imath}$ (for the case of the pericenter just replace $e$ by $-e$). In our case $G=1$, $m = 3$,  $m_4 = 0$, therefore
$$
\begin{array}{c}
{\bf r} = {\bf 0}, \enskip {\bf v} = {\bf 0}, \\[6pt]
{\bf r}_4 = 
\left(
\begin{array}{c}
a(1+e) \\
0
\end{array}
\right),
\enskip
{\bf v}_4 = 
\left(
\begin{array}{c}
0 \\
\displaystyle  \sqrt{\frac{3(1-e)}{a(1+e)}}
\end{array}
\right),
\end{array}
$$
so that
\begin{equation}
v_{y4} = \frac{\sqrt{3(1-e)}}{\sqrt{x_4}}.
\label{aprox}
\end{equation}
The curve defined by (\ref{aprox}) resembles the first graph of Fig. \ref{fig5}. The orbits defined by the two-body problem are periodic. However, in order to fulfill the relation between the orbits in the restricted four-body problem, we only need to consider those whose period is an even multiple of the characteristic time, that is $T = 2m\overline{T}$, $m \in \mathbb{N}$. With this we obtain the relation

\begin{equation}
    ( 2m\overline{T})^2 = \frac{4}{3} \pi^2 \frac{x_4^3}{(1+e)^3}.
    \label{restriccionperiodoE}
\end{equation} 
Notice that (\ref{aprox}) is invariant under the transformation
$$
(x_4,v_{y4}) \to (\alpha x_4,\frac{1}{\sqrt{\alpha}} v_{y4}),
$$
which corresponds to a change of the scale of the orbit (homothetic property). It implies that, for a given eccentricity $e$, the entries $x_4$ and $v_{y4}$ are related by a monoparametric family. However, in our case, the scale is fixed because the eight figure orbit has specific initial conditions. Thus, the periodic orbits do not appear in a monoparametric family, only in a discrete way. It is expected that the curve (\ref{aprox}), for a large distance from the test particle to the primaries, gives rise to good enough approximations of initial conditions of periodic orbits of the restricted four-body problem.

\section{Conclusion}
\label{VII}

We have introduced a restricted four-body problem, associated to the eight figure choreography, and we showed the existence of a reversing symmetry for the $N$--body problem with $N=2n+k$, where each pair of the $2n$ bodies have equal masses. With the help of reversing symmetries, we computed a small part of the sets that give rise to symmetric orbits.

We approximate the solutions with a two-body problem, when the test particle is far away from the primaries. In principle the same idea could work for orbits in which the test particle surround only one primary. Solutions which pass near collision are interesting, in particular those associated to binary systems. Nevertheless, their numerical study require from the regularization of binary collisions since the test particle passes very close to some of the choreographic bodies. In the future, this work can be complemented with a more extensive study of symmetric periodic orbits, and its bifurcations. 

Another kind of continuation that can be made in this restricted four-body problem is in the mass parameter $m_4$, with the aim of connecting orbits of the full gravitational problems of three and four bodies. For instance, it would be interesting to relate the eight figure and Gerver's super eight choreographies.

\section*{ACKNOWLEDGMENTS}
The first author is pleased to acknowledge the financial support of Instituto Tecnol\'ogico Aut\'onomo de M\'exico (ITAM). The second author is pleased to acknowledge the support of Asociaci\'on Mexicana de Cultura A.C., and National System of Researchers (SNI).

\endpaper

\end{document}